\documentclass[a4paper,11pt, leqno]{article}
\usepackage{amsthm}
\usepackage[francais,english]{babel}
\usepackage[intlimits]{amsmath}
\usepackage{amssymb}
\usepackage{latexsym}
\usepackage{mathrsfs}
\usepackage[all]{xy}
\RequirePackage{ifthen}
\RequirePackage{verbatim}

\newcounter{segno}[section]

\renewcommand{\thesegno}{\thesection.\arabic{segno}}

\newskip\segskipamount
\newskip\procskipamount
\newskip{\interskipamount}
\newcommand{\segskip}{\vskip\segskipamount}
\newcommand{\procskip}{\vskip\procskipamount}
\newcommand{\interskip}{\vskip\interskipamount}

\newcommand{\segbreak}{\par
   \ifdim\lastskip<\segskipamount\removelastskip
   \penalty-200
   \segskip\fi}

\newcommand{\procbreak}{\par
   \ifdim\lastskip<\procskipamount\removelastskip
   \penalty-100
   \procskip\fi}

\newenvironment{segment}{%
  \segbreak
  \refstepcounter{segno}%
  \noindent\textbf{(\thesegno)\ }}%
  {\par}

\newenvironment{segproc}[1]{%
  \textbf{#1.\ }%
  \em}%
  {\par\interskip}

\newenvironment{proclamation}[1]{%
  \procbreak
  \noindent\textbf{#1.\ }%
  \em}%
  {\par\interskip}

  {\par}%

\numberwithin{equation}{segno}

\newcounter{listcounter}
\newcounter{deflistcounter}

\newskip{\itemsepamount}
\newskip{\topsepamount}
\newenvironment{assertionlist}{%
  \begin{list}
    {\upshape (\arabic{listcounter})}
    {\setlength{\leftmargin}{18pt}
     \setlength{\rightmargin}{0pt}
     \setlength{\itemindent}{0pt}
     \setlength{\labelsep}{5pt}
     \setlength{\labelwidth}{13pt}
     \setlength{\listparindent}{\parindent}
     \setlength{\parsep}{0pt}
     \setlength{\itemsep}{\itemsepamount}
     \setlength{\topsep}{\topsepamount}
     \usecounter{listcounter}}}
  {\end{list}}

\newcommand{\FF}{\mathbb{F}}

\newcommand{\ZZ}{\mathbb{Z}}

\newcommand{\wbar}{\bar{w}}

\newcommand{\Gbar}{\bar{G}}

\newcommand{\pgbar}{\bar{\pi}}

\newcommand{\wdot}{\dot{w}}

\newcommand{\htilde}{\tilde{h}}

\newcommand{\Htilde}{\tilde{H}}

\newcommand{\FFbar}{\overline{\FF}}

\newcommand{\gfr}{{\mathfrak g}}

\newcommand{\Hscr}{{\mathscr H}}

\newcommand{\Acal}{{\mathcal A}}

\newcommand{\Ocal}{{\mathcal O}}

\newcommand{\Xcal}{{\mathcal X}}

\newcommand{\eps}{\varepsilon}

\newcommand{\dgtilde}{\tilde\delta}

\DeclareMathOperator{\DR}{DR}

\DeclareMathOperator{\Fr}{Fr}

\DeclareMathOperator{\GL}{GL}

\DeclareMathOperator{\Ker}{Ker}
\DeclareMathOperator{\Lie}{Lie}

\DeclareMathOperator{\opp}{opp}

\DeclareMathOperator{\relpos}{relpos}
\DeclareMathOperator{\Ref}{Ref}

\DeclareMathOperator{\Sp}{Sp}
\DeclareMathOperator{\Spec}{Spec}

\newcommand{\ad}{\textnormal{ad}}

\newcommand{\lrangle}{{\langle\ ,\ \rangle}}

\newcommand{\restricted}[1]{{}_{\vert}{}_{#1}}
\newcommand{\set}[2]{\{\,#1 \mid #2\,\}}

\newcommand{\twosmallmatrix}[4]{\left(\begin{smallmatrix} #1 & #2 \\ #3 &
      #4 \end{smallmatrix}\right)}

\newcommand{\lto}{\longrightarrow}
\newcommand{\bijective}{\leftrightarrow}
\newcommand{\epi}{\twoheadrightarrow}

\newcommand{\sends}{\mapsto}

\newcommand{\iso}{\overset{\sim}{\to}}
\newcommand{\liso}{\overset{\sim}{\lto}}

\begin{document}

\title{Specialization of $F$-Zips}

\author{Torsten Wedhorn}

% \address{Mathematisches Institut Bonn\\
% Beringstra\ss e 4\\
% 53115 Bonn\\
% Germany}

% \email{wedhorn@math.uni-bonn.de}

% \urladdr{http://www.math.uni-bonn.de/people/wedhorn}

%\date{August 15, 2004}

\maketitle

\noindent{\scshape Abstract.\ }
In \cite{MW}, B.~Moonen and the author
defined a new invariant, called $F$-Zips, of certain varieties in positive
characteristics. We showed that the isomorphism
classes of these invariants can be interpreted as orbits of a certain
variety $Z$ with an action of a reductive group $G$. In
loc.~cit.\ we gave a combinatorial description of the set of these orbits.

In this manuscript we give an explicit combinatorial recipe to decide
which orbits are in the closure of a given orbit. We do this by
relating $Z$ to a semi-linear variant of the wonderful
compactification of $G$ constructed by de Concini and Procesi. As an
application we give an explicit criterion of the closure relation for
Ekedahl-Oort strata in the moduli space of principally polarized
abelian varieties.

%==========================================================================

\section*{Introduction}

Let $k$ be a field of characteristic $p > 0$ and let $X$ be a smooth,
proper scheme over $k$ such that the Hodge spectral sequence
degenerates. Examples are abelian varieties, complete
intersections in a projective space, K3-surfaces, toric varieties, and
curves. In \cite{MW} it was shown
that for each $i$ the $i$-th de Rham cohomology of $X$ is endowed with the
structure of an $F$-zip (see loc.~cit., Definition~2.1, for the
definition of an $F$-zip). Moreover, $F$-zips also arise as mod $p$
reductions of completely divisible lattices in filtered isocrystals
(see loc.~cit., Example~7.3) and from the $p$-torsion of a
Barsotti-Tate group. The main result in loc.~cit.\ then is the
classification of isomorphism classes of $F$-zips in terms of
combinatorial data associated to certain Coxeter groups.

This result was obtained by showing that isomorphism classes of
$F$-zips correspond to $G$-orbits of a certain variety $Z_{G,J}$. Here
$G = \GL_n$ for some $n$ and $J$ is a fixed set of simple reflections
in the Weyl group of $G$ (see loc.~cit., Section~4.8, or \eqref{defZ}
below for the definition of $Z_{G,J}$). In fact, the varieties
$Z_{G,J}$ can be defined for an arbitrary reductive group $G$, defined
over a finite field, and for any set $J$ of simple reflections in the
Weyl group $W$ of $G$. On $W$ there is a natural partial order, namely
the Bruhat order, which we denote by $\leq$. Then in loc.~cit.\ it was
shown that there is a natural bijection
\[
\{\text{$G$-orbits of $Z_{G,J}$}\} \bijective {}^JW,
\]
where ${}^JW = \set{w \in W}{\text{$w < sw$ for all $s \in
    J$}}$. In particular, there are only finitely many $G$-orbits in $Z_{G,J}$.

In this paper we will prove how these orbits specialize into each
other: For $w \in {}^JW$ denote the corresponding $G$-orbit of
$Z_{G,J}$ by $Z_{G,J}^w$. Let $x = w_0^J$ be the maximal element in
$W^J = \set{w \in W}{\text{$w < ws$ for all $s \in J$}}$, and let
$W_J$ be the subgroup of $W$ generated by $J$. We denote by $F\colon W
\iso W$ the automorphism of $W$ induced by the Frobenius. Then the
main result is (see \eqref{closureZ}):

\begin{proclamation}{Theorem}
For any $w,w' \in {}^JW$ the following two assertions are equivalent:
\begin{assertionlist}
\item $Z_{G,J}^w \subset \overline{Z_{G,J}^{w'}}$.
\item There exists $u \in W_J$ such that $u^{-1}w\delta(u) \leq w'$
  with $\delta(u) = xF(u)x^{-1}$.
\end{assertionlist}
\end{proclamation}

We obtain the following application to the Ekedahl-Oort stratification
of the moduli space $\Acal_g$ of $g$-dimensional principally polarized
abelian varieties in characteristic $p$ (or more generally to good
reductions of Shimura varieties of PEL-type): One corollary of the
results in \cite{MW} is that the Ekedahl-Oort strata of $\Acal_g$
(e.g., defined in \cite{Oort}) can be parametrized by ${}^JW$. Here
$W$ is the Weyl group of the symplectic group $\Sp_{2g}$ and
$J$ is the type of the Siegel parabolic in $\Sp_{2g}$, see
\eqref{symplectic} for the precise definition of $W$ and $J$. This
parametrization had been obtained beforehand by F.~Oort in \cite{Oort}
with a different formulation and, for arbitrary good
reductions of Shimura varieties of PEL-type, by B.~Moonen in
\cite{Mo}. Moreover, F.~Oort (\cite{Oort}) and the author
(\cite{WdOortstrata}) both have shown by different methods that the
Ekedahl-Oort strata are locally closed and equidimensional and that
the closure of an Ekedahl-Oort stratum is a union of Ekedahl-Oort
strata. But it remained an open question which strata are contained in
the closure of a given one (see \cite{Open}, Problem~11). Now the
theorem above implies:

\begin{proclamation}{Corollary}
For $w \in {}^JW$ denote by $\Acal_g^w$ the corresponding Ekedahl-Oort
stratum in $\Acal_g$. For $w,w' \in {}^JW$ the following two
assertions are equivalent:
\begin{assertionlist}
\item $\Acal_g^w \subset \overline{\Acal_g^{w'}}$.
\item There exists $u \in W_J$ such that $u^{-1}w\delta(u) \leq w'$
  with $\delta(u) = w_{0,J}uw_{0,J}$.
\end{assertionlist}
Here $w_{0,J}$ denotes the maximal element in $W_J$.
\end{proclamation}

The varieties $Z_{G,J}$ are a Frobenius-linear variant of varieties $Z'_{G,J}$
defined by Lusztig in \cite{Lu1} and \cite{Lu2}. In this linear
setting he also obtains $G$-stable subvarieties $Z^{\prime w}_{G,J}$
which are parametrized (after renormalization) by ${}^JW$. Also, X.~He
shows in \cite{He}, that the closure of $Z^{\prime w}_{G,J}$ is a
union of such $G$-stable subvarieties. Although the
constructions in \cite{Lu1}, \cite{Lu2}, and \cite{He} do not carry
over to the Frobenius-linear setting, the proof of the Theorem above
uses similar methods.

We introduce notations in Section~1 and recall some
definitions and results from \cite{MW} in Section~2. Then it is shown
that $Z_{G,J}$ carries a transitive $(G \times G)$-action. Moreover
there is a $(G \times G)$-equivariant smooth surjective morphism from
$Z_{G,J}$ to a Frobenius-linear variant of a $(G \times G)$-orbit
$X_{G^{\rm ad},J}$ in the wonderful compactification by de Concini and
Procesi of the adjoint group $G^{\rm ad}$. This morphism induces
a bijection of $(B \times B)$-orbits, where $B \subset G$ is a Borel
subgroup. As the closure relation of the $(B \times B)$-orbits in the
wonderful compactification is
known by the work of Springer \cite{Sp}, one obtains the
closure relation of the $(B \times B)$-orbits in $Z_{G,J}$. Moreover,
for every $w \in {}^JW$ we find a $(B \times B)$-orbit which is
contained in $Z^w_{G,J}$, see Section~3.

In Section~4 the partial order on ${}^JW$ is introduced which
corresponds to taking closures of the $G$-orbits $Z^w_{G,J}$. Although we have
to work in a more general setting than He in \cite{He}, the proofs in
this section are easy modifications of the proofs given by He. Finally
in Section~5, we prove the main theorem making use of the connection
of $G$-orbits and $(B \times B)$-orbits.

In the last section we apply these results to the Ekedahl-Oort
stratification and we obtain the corollary above.

%--------------------------------------------------------------------------

\section{Refinement of parabolic subgroups}

\begin{segment}\label{defp}
In the sequel $p$ is a prime
number and $q$ is a fixed power of~$p$. For a scheme~$S$ of
characteristic~$p$ we denote by $\Fr_S\colon S \to S$ the morphism
which is the identity on the underlying topological space and the
homomorphism $x \mapsto x^q$ on the sheaves of rings. For an
$\Ocal_S$-module $M$ we set $M^{(q)} = \Fr_S^*M$.

Let $\FFbar$ be an algebraic closure
of~$\FF_q$. All varieties we will consider are reduced schemes of
finite type over $\FFbar$, although sometimes they will have fixed
$\FF_q$-rational structures. We will also systematically confuse a
variety $X$ and its $\FFbar$-valued points.
\end{segment}

\begin{segment}\label{defgroup}
Let $G$ be a connected reductive group over~$\FF_q$. We fix a maximal
torus $T$ and a Borel subgroup $B$ of $G$ containing $T$ which are both
defined over $\FF_q$. We denote by $G^{\ad}$ the adjoint group of $G$
and by $T^{\ad}$ and $B^{\ad}$ the image of $T$ and $B$ in
$G^{\ad}$. More general, we write $H^{\ad}$ for the image in $G^{\ad}$
of an algebraic subgroup $H$ of $G$.

Let $W$ be the Weyl group of $G$ associated to
$T$ and let $I \subset W$ be the set of simple reflections
corresponding to $B$. We endow $W$ and any subset of $W$ with the
Bruhat order which will be denoted by $\leq$ and denote by $w \sends
\ell(w)$ the length function on $W$.

If there is no risk of confusion, we simply write $F\colon G
\to G$ for the Frobenius $F_{G/\FF_q}\colon G \to
G^{(q)} = G$, which is the map $x \sends x^q$ on local sections.
It is an endomorphism of~$G$ which induces an
automorphism of the Weyl group $W$, again denoted by~$F$. As $B$ is
defined over $\FF_q$, we have $F(I) = I$.

For any subset $J \subset I$ we denote by $W_J$ the subgroup of $W$
generated by $J$. We denote by $w_{0,J}$ the longest element in $W_J$
and set $w_0 = w_{0,I}$. Finally let $w_0^J = w_0w_{0,J}$ be the
element of maximal length in $W^J$. Similarly we set ${}^Jw_0 =
w_{0,J}w_0$. Finally we write $J^{\opp} = w_0Jw_0 \subset I$.

We denote by $P_J \supset B$ the standard parabolic subgroup
corresponding to $J$ and by $L_J \subset P_J$ the unique
Levi subgroup of $P_J$ with $L_J \supset T$. A parabolic subgroup $P$
of $G$ is called of \emph{type $J$} if $P$ is conjugated to $P_J$.

Let $P_J^{\opp}$ be the unique parabolic subgroup of $G$ such that
$P_J \cap P_J^{\opp} = L_J$. It is a parabolic subgroup of type
$J^{\opp}$ such that $w_0P_J^{\opp}w_0 = P_{J^{\opp}}$. We denote by
$\pi_J$ and $\pi_J^{\opp}$ the projections of $P_J$ (resp.~$P_J^{\opp}$) onto
$L_J$.

For any parabolic subgroup $P$ of $G$ we denote by $U_P$ its unipotent
radical. We set $U = U_B$ and $U^{\rm opp} = U_{B^{\opp}}$. For each subset $J$
of $I$ we set $U_J = L_J \cap U$ and $U_J^{\opp} = L_J \cap
U^{\opp}$. Then multiplication induces isomorphisms
\begin{equation}\label{multunipotent}
U_J \times U_{P_J} \liso U \overset{\sim}{\longleftarrow} U_{P_J}
\times U_J.
\end{equation}

Let $\Phi \subset X^*(T)$ be the set of roots of $(G \otimes_{\FF_q}
\FFbar,T \otimes_{\FF_q} \FFbar)$. For each $i \in I$ let $\alpha_i$
be the corresponding simple root in $\Phi$. We denote by $\Phi_J$
the set of roots which are in the $\ZZ$-span of $\set{\alpha_j}{j
\in J}$ and set $\Phi_J^+ = \Phi^+ \cap \Phi_J$ where $\Phi^+$ is
the set of positive roots.
\end{segment}

\begin{segment}\label{defineJW}
We denote by ${}^JW$ (resp.~$W^J$) the set of minimal length coset
representatives of $W_J\backslash W$ (resp.~$W/W_J$). For $J,K \subset
I$ we write ${}^JW^K = {}^JW \cap W^K$. This is a system of
representatives in $W$ for $W_J\backslash W/W_K$.

We have the following descriptions for ${}^JW$:
\begin{equation}\label{describeJW}
\begin{aligned}
{}^JW &= \set{w \in W}{\text{$w \leq w'$ for all $w' \in
    W_Jw$}}\\
&= \set{w \in W}{\text{$\ell(w) < \ell(sw)$ for all $s \in J$}}\\
&= \set{w \in W}{\text{$\ell(uw) = \ell(u) + \ell(w)$ for all $u \in
    W_J$}}\\
&= \set{w \in W}{w^{-1}(\Phi_J^+) \subset \Phi^+}.
\end{aligned}
\end{equation}
For every element $x \in W$ there exist unique elements $u \in W_J$
and $w \in {}^JW$ such that $x = uw$. We call $u$ the \emph{$W_J$-part
  of $x$} and $w$ the \emph{${}^JW$-part of $x$}.

If $J' \subset J \subset I$ are two subsets, we clearly have ${}^JW
\subset {}^{J'}W$. Conversely, we also have a canonical surjective map
\begin{equation}\label{canonicalWeyl}
{}^{J'}W \bijective W_{J'}\backslash W \epi W_J\backslash W \bijective {}^JW.
\end{equation}
\end{segment}

\begin{segment}\label{JWremark}
Let $w \in {}^JW$ and $s \in I$. By \cite{Bourbaki}, \S 1,
$\textnormal{n}^{\textnormal{o}}$ 1.7, there are three possibilities:
\begin{assertionlist}
\item $ws > w$ and $ws \in {}^JW$.
\item $ws > w$ and $ws = tw$ for some $t \in J$.
\item $ws < w$ and $ws \in {}^JW$.
\end{assertionlist}
\end{segment}

\begin{segment}
We have the following easy lemma.

\begin{proclamation}{Lemma}\label{convertlemma}
Let $J \subset I$ and set $K = J^{\opp} = w_0Jw_0$.
\begin{assertionlist}
\item
The map $x \mapsto x^{-1}$ defines a bijection ${}^JW \bijective W^J$ which
preserves Bruhat order and length.
\item
The map $x \mapsto w_0xw_0$ defines a bijection $W^J \bijective W^K$ which
preserves Bruhat order and length.
\item
The map $x \mapsto x^{-1}w_0^K$ defines a bijection ${}^JW \to W^K$
which reverses Bruhat order. Moreover $\ell(x^{-1}w_0^K) =
\ell(w_0^K) - \ell(x)$.
\end{assertionlist}
\end{proclamation}
\end{segment}

\begin{segment}\label{defrelpos}
For any subgroup $H$ of $G$ and for any element $g \in G$ we set
${}^gH = gHg^{-1}$.

For any two parabolic subgroups $P$ and $Q$ of type $J$ and $K$,
respectively, we denote by $\relpos(P,Q) \in {}^JW^K \subset W$
their relative position. Recall that the relative position of $P$ and
$Q$ can be defined as follows: Denote by $W(P,Q)$ the set of elements
$w \in W$ such that there exists a $g \in G$ such that ${}^gP = P_J$
and ${}^gQ = {}^wP_K$. Then $W_J$ acts from the left and $W_K$ acts
from the right on $W(P,Q)$, and $W(P,Q)$ is a single orbit for this
action. Hence there exists a unique element $\relpos(P,Q) \in W(P,Q)
\cap {}^JW^K$.

In particular we have
\begin{alignat}{2}
\relpos({}^gP,{}^gQ) &= \relpos(P,Q) && \qquad\text{for all $g \in G$,}\\
\relpos(P,Q) &= \relpos(Q,P)^{-1},\\
\relpos(P_J^{\opp},P_J) &= w_0^J,\\
\relpos(P_J,{}^wP_K) &= w && \qquad\text{for all $w \in {}^JW^K$.}
\end{alignat}
\end{segment}

\begin{segment}\label{refinepsgp}
If $P$ and $Q$ are two parabolic subgroups of $G$, we set
\[
\Ref_Q(P) = (P \cap Q)U_P,
\]
the refinement of $P$ by $Q$ (cf.~\cite{MW}, Section~3.7 and
Example~3.3). This is a parabolic subgroup of $G$ which is contained in $P$.

If $J$ and $K$ are the types of $P$ and $Q$, respectively, $\Ref_Q(P)$
is of type $J \cap {}^{\relpos(P,Q)}K$.
\end{segment}

\begin{segment}\label{refinelemma}
\begin{segproc}{Lemma}
Let $J$ and $K$ be two subsets of $I$, and let $w \in W$. Then
$\Ref_{{}^wP_K}(P_J)$ contains $B$ if and only if
\[
{}^{w^{-1}}\Phi_J^+ \subset \Phi_K \cup \Phi^+.
\]
In particular, this is the case if $w \in {}^JW$.
\end{segproc}

\begin{proof}
For any root $\alpha$ let $\gfr_{\alpha} \subset \Lie(G)$ be the
subspace where $T$ acts via $\alpha$. As $P := \Ref_{{}^wP_K}(P_J)$
contains the chosen maximal torus $T$, we can write $\Lie(P) = \Lie(T)
\oplus \bigoplus_{\alpha \in \Phi_P} \gfr_{\alpha}$ for a certain
subset $\Phi_P$ of $\Phi$. Then $B$ is contained in $P$ if and only if
$\Phi^+ \subset \Phi_P$. By definition of $P$ we have
\[
\Phi_P = \bigl((\Phi_J \cup \Phi^+) \cap ({}^w\Phi_K \cup
{}^w\Phi^+)\bigr) \cup (\Phi^+ \setminus \Phi_J)
\]
and therefore $\Phi^+ \subset \Phi_P$ if and only if $\Phi_J^+ \subset
{}^w\Phi_K \cup {}^w\Phi^+$. The last claim follows from \eqref{describeJW}.
\end{proof}
\end{segment}

\begin{segment}\label{Howlett}
We have the following result by Howlett (see \cite{Carter},
Proposition 2.7.5):

\begin{proclamation}{Lemma}
Let $J,K \subset I$ and $\wbar \in {}^JW^K$. Set $K' := K \cap
  {}^{\wbar^{-1}}J$. Then each element $w \in
W_J\wbar W_K$ can be uniquely expressed in the form $w = u\wbar v$
with $u \in W_J$ and $v \in W_K \cap {}^{K'}W$. Moreover, we have
  $\wbar v \in {}^JW$ and
  $\ell(w) = \ell(u) + \ell(\wbar) + \ell(v)$.
\end{proclamation}
\end{segment}

\begin{segment}\label{Howlettvariant}
In a special case, we have also the following more precise version of
Howlett's lemma, proved in \cite{He}, Lemma~3.6:

\begin{proclamation}{Lemma}
Let $J,K \subset I$ and $\wbar \in {}^JW^K$. Set $K' := K \cap
  {}^{\wbar^{-1}}J$ and $J' := J \cap {}^{\wbar}K$. Let $w \in \wbar W_K$. Then
the unique decomposition of $w$ in Howlett's lemma is of the form
$w = u\wbar v$ with $v \in W_K \cap {}^{K'}W$
and $u \in W_{J'}$.
\end{proclamation}
\end{segment}

%--------------------------------------------------------------------------

\section{The classifying variety}

\begin{segment}\label{defJ}
From now on we fix a subset $J \subset I$ which is defined over
$\FF_q$, in particular we have $F(P_J) = P_J$ and $F(J) = J$. We set
$K = w_0Jw_0 \subset I$ which is also defined over $\FF_q$.
\end{segment}

\begin{segment}\label{defZ}
Let $Z_J = Z_{G,J}$ be the variety of triples $(P,Q,[g])$ where $P$ and $Q$ are
parabolic subgroups of types $J$ and $K$, respectively, and where
$[g]$ is a double coset in $U_Q\backslash G/F(U_P)$ such that
$\relpos(Q,{}^gF(P)) = w_0^J = {}^Kw_0$. This is a special case of the
variety $Z_J$ defined in \cite{MW}, Section~4.8.

The variety $Z_J$ carries a left $G \times G$-action by
\[
(h,h')\cdot (P,Q,[g]) = ({}^{h'}P,{}^hQ,[hgF(h')^{-1}]).
\]
\end{segment}

\begin{segment}\label{definesequence}
For any point $(P,Q,[g]) \in Z_J$ we define a sequence of pairs of
parabolic subgroups $(P_n,Q_n)_{n\geq 0}$ inductively:
\begin{equation}\label{seqparabolic}
\begin{aligned}
P_0 &:= P, & \quad Q_0 &:= Q, \\
P_n &:= \Ref_{Q_{n-1}}(P_{n-1}), & \quad Q_n &:=
\Ref_{{}^gF(P_n)}(Q_{n-1}).\\
\end{aligned}
\end{equation}
Clearly $P_n \subset P_{n-1}$ and $Q_n \subset Q_{n-1}$. Therefore
there exists an $N \geq 0$ such that $(P_n,Q_n) = (P_{n+1},Q_{n+1})$
for all $n \geq N$. We set $(P_\infty,Q_\infty) = (P_N,Q_N)$ and obtain
\[
\relpos(P_\infty,Q_\infty) \in {}^{J_\infty}W^{K_\infty}
\]
where $J_\infty \subset J$ and $K_\infty \subset K$ are the types of
$P_\infty$ and $Q_\infty$, respectively. It follows from \cite{MW},
Theorem~4.11 and Section~4.6, that
$\relpos(P_\infty,Q_\infty) \in {}^JW$ and we set $\alpha([g,g',z]) :=
\relpos(P_\infty,Q_\infty)$. In this way we obtain a map
\[
\alpha\colon Z_J \lto {}^JW.
\]
\end{segment}

\begin{segment}\label{GorbitsofZ}
Now Theorem~4.11 of \cite{MW} describes the $G$-orbits
of $Z_J$ where we consider $G$ embedded
diagonally in $G \times G$. We get:

\begin{proclamation}{Theorem}
The map $\alpha$ induces a bijection of the set of $G$-orbits of $Z_J$
and the set ${}^JW$.
\end{proclamation}

We denote the $G$-orbit of $Z_J$ corresponding to $w \in
{}^JW$ by $Z_J^w$.
\end{segment}

\begin{segment}\label{Jinfty}
Let $z \in Z_J$ and let $(P_n,Q_n)_n$ be the associated
sequence of pairs of parabolic subgroups \eqref{GorbitsofZ} and let
$J_{\infty}$ (resp.~$K_{\infty}$) be the type of $P_{\infty}$
(resp.~$Q_{\infty}$). Assume that $z \in Z_J^w$ for some $w \in
{}^JW$. It follows from \cite{MW}, Section~4.6, that
$J_{\infty}$ can be described as the largest subset $J'$ of $J$ such
that ${}^{ww_0^J}J' = J'$ and that $K_{\infty} = {}^{w^{-1}}J_{\infty}
= {}^{w_0^J}J_{\infty}$.
\end{segment}

\begin{segment}\label{keylemma}
For any $w \in W$ we choose an element $\wdot \in N_G(T)$ which maps
to $w$. We set
\begin{equation}\label{definedelta}
\delta\colon W_J \liso W_K, \qquad u \sends w_0^JF(u)(w_0^J)^{-1}.
\end{equation}

\begin{proclamation}{Lemma}
Let $w \in W$ and $b \in B$ and set
\[
z = z(w,b) = (P_J,{}^wP_K,[\wdot \wdot^J_0 F(b)]) \in Z_J.
\]
We write $w = uw'$ with $u \in W_J$ and $w' \in {}^JW$. Then there
exists an element $v \in W_K$ such that $w'v \in
{}^JW$, $z \in Z^{w'v}_J$ and $\delta^{-1}(v) \leq u$.
\end{proclamation}
\end{segment}

\begin{segment}
\begin{segproc}{Corollary}\label{keycorollary}
Let $w \in {}^JW$ and $b \in B$ and set
\[
z = z(w,b) = (P_J,{}^wP_K,[\wdot \wdot^J_0 F(b)]) \in Z_J.
\]
Then $z \in Z^w_J$.
\end{segproc}

\begin{proof}
With the notations of \eqref{keylemma}, we have $u = 1$ and hence
$v=1$.
\end{proof}
\end{segment}

\begin{segment} {\it Proof of \eqref{keylemma}}.\ 
Let $(P_n,Q_n)_{n \geq 0}$ be the sequence of pairs of parabolics
associated to $z$ \eqref{definesequence}. Let $J_n$ be the type of
$P_n$ and $K_n$ be the type of $Q_n$. We set $y_n =
\relpos(P_n,Q_n)$. Let ${\infty} \geq 0$ be the smallest integer such that
$P_{\infty} = P_{{\infty}+1}$. Then we have for all $n \geq {\infty}$,
$P_n = P_{\infty}$, $Q_n = Q_{\infty}$, $J_n = J_{\infty}$, $K_n =
K_{\infty}$ and $y_n = y_{\infty}$.

By definition $z \in Z_J^{y_{{\infty}}}$ and
$y_{{\infty}} \in {}^JW^{K_{{\infty}}}$. We have
to show that $y_{{\infty}} = w'v$ for some $v \in W_K$ with
$\delta^{-1}(v) \leq u$. For all $n \geq 0$,
\[
J_{n+1} = J_n \cap {}^{y_n}K_n, \qquad K_{n+1} = {}^{w_0^J}J_{n+1}
\]
and $K_{{\infty}} = {}^{w_0^J}J_{{\infty}} =
{}^{y_{{\infty}}^{-1}}J_{{\infty}}$ by \cite{MW}, Section~4.6.

We first consider two special cases. The first one is that ${\infty} = 0$.
Then $y_{\infty} = \relpos(P,Q)$. By \eqref{Howlett} we can  write $w' =
y_{\infty}v$ with $v \in W_K \cap {}^{K \cap
  {}^{y_{\infty}^{-1}}J}W$. As $J_1 = J_0$, we have $K \cap
{}^{y_{\infty}^{-1}}J = K$ and therefore $v = 1$. This
proves the lemma in this case.

The second special case is the case $u = 1$, i.e., we prove the
Corollary~\eqref{keycorollary}. Then we have
to show that $z \in Z^w_J$. We claim
that $B \subset P_n \cap {}^{w^{-1}}Q_n$ for all $n \geq 0$, in other
words, $P_n = P_{J_n}$ and $Q_n = {}^wP_{K_n}$. The claim is shown by
induction on $n$.

It certainly holds for $P_0 = P_J$ and $Q_0 = {}^wP_K$. Assume that it
holds for $n \geq 0$. Then we have
\[
P_{n+1} = \Ref_{{}^wP_{K_n}}(P_{J_n}) \supset B
\]
by \eqref{refinelemma}, as $w \in {}^JW \subset {}^{J_n}W$. Moreover
\[
{}^{w^{-1}}Q_{n+1} = {}^{w^{-1}}\Ref_{{}^{w w^J_0
    b}F(P_{J_{n+1}})}(Q_n) = \Ref_{{}^{w^J_0}F(P_{J_{n+1}})}(P_{K_n})
    \supset B
\]
again by \eqref{refinelemma}, as $w_0^J \in {}^KW^J \subset
{}^{K_n}W$. This proves the claim.

Therefore we have seen that $w \in
W_{J_{{\infty}}}y_{{\infty}}W_{K_{{\infty}}} \cap {}^{J_{{\infty}}}W$. By
\eqref{Howlett} this implies $w = y_{{\infty}}v$ for some $v \in
W_{K_{{\infty}}} \cap {}^{K_{\infty} \cap
  {}^{y_{\infty}^{-1}}J_{\infty}}W = W_{K_{{\infty}}} \cap
{}^{K_{{\infty}}}W$. Hence $v = 1$ and $w = y_{{\infty}}$, which
proves the case $u = 1$.

In general we define $Z_{J_n,w^0_J}$ as the variety of triples
$(P,Q,[g])$, where $P$ is a parabolic subgroup of type $J_n$, $Q$ is a
parabolic subgroup of type $K_n = {}^{w^J_0}J_n$, and $[g] \in
U_Q\backslash G/U_{F(P)}$ with $\relpos(Q,{}^gF(P)) = w_0^J$. This
variety carries a $G \times G$-action as in \eqref{defZ}, and in
particular a $G$-action where we consider $G$ embedded diagonally in $G
\times G$. As in \eqref{GorbitsofZ}, the $G$-orbits of $Z_{J_n,w^0_J}$
are parametrized by ${}^{J_n}W$ (\cite{MW}, Theorem~4.11). In particular,
we can consider the $G$-orbit $Z^{y_{{\infty}}}_{J_n,w_0^J}$
corresponding to the element $y_{{\infty}} \in {}^JW \subset {}^{J_n}W$. By
\cite{MW}, Lemma~4.6 and Lemma~4.7, for any $n \geq 0$ the map
\begin{align*}
\vartheta\colon Z^{y_{{\infty}}}_{J_n,w_0^J} &\to
Z^{y_{{\infty}}}_{J_{n+1},w_0^J},\\
(P,Q,[g]) &\sends (\Ref_Q(P),\Ref_{{}^gF(\Ref_Q(P))}(Q),[g])
\end{align*}
is a well defined surjective smooth morphism which induces an
isomorphism of the fppf-quotients
\[
G\backslash Z^{y_{{\infty}}}_{J_n,w_0^J} \liso G\backslash
Z^{y_{{\infty}}}_{J_{n+1},w_0^J}.
\]

% By \cite{MW}, Section~4.6, there exists an element $v_n \in W_{K_n}$ such that
% $y_{n+1} = y_nv_n$. As $y_n \in {}^{J_n}W^{K_n}$, it follows from
% \eqref{Howlettvariant} that $v_n \in W_{K_n} \cap {}^{K_{n+1}}W$.

We define inductively elements $u_n \in W_{J_n} \cap W^{J_{n+1}}$ and
$u'_n \in W_{J_n}$ as follows: For $n = 0$ we set $u'_0 := u$ and
write $u = u'_0 = u_0u'_1$ with $u_0 \in W_J \cap W^{J_1}$ and $u'_1
\in W_{J_1}$. For $n \geq 1$ we write $u'_n = u_nu'_{n+1}$ with $u_n
\in W_{J_n} \cap W^{J_{n+1}}$ and $u'_{n+1} \in W_{J_{n+1}}$. This
completes the inductive definition. Then $u_n = 1$ for $n \geq {\infty}$ and
\begin{equation}\label{addu}
u = u_0u_1u_2\cdots u_{{\infty}-1}, \qquad \ell(u) = \ell(u_0) + \ell(u_1) +
\cdots + \ell(u_{\infty-1}).
\end{equation}

For the proof of the lemma we can replace $z$ by
\[
z_0 = (P_J,{}^{w'}P_K,[w'w_0^J F(b) F(u)]),
\]
as both are in the same $G$-orbit. In general, let $z_n \in
Z_{J_n,w_0^J}$ be an element of the form
\[
z_n = (P_{J_n},{}^{w'_n}P_{K_n},[h_n w'_nw_0^J F(b_n) F(u'_n)])
\]
where $h_n \in {}^{w'}U_K$, $w'_n \in {}^{J_n}W \cap w'W_K$, $b_n \in
B$. By definition, we can write $u'_n =
u_nu'_{n+1}$. By \eqref{commutelemma} below there exist $v_n \in
W_K$, $h_{n+1} \in {}^{w'}U_K$, and $b_{n+1} \in B$ such that
$\delta^{-1}(v_n) \leq u_n$, $w'_{n+1} := w'_nv_n \in {}^{J_{n+1}}W$,
and such that $\vartheta(z_n)$ is in the same $G$-orbit as
\[
z_{n+1} := (P_{J_{n+1}},{}^{w'_{n+1}}P_{K_{n+1}},[h_{n+1} w'_{n+1}
w_0^J F(b_{n+1}) F(u'_{n+1})]).
\]

If we start this process with $z_0$ and $w'_0 = w'$, we get a sequence
of elements $w'_0,w'_1,\dots$ with $w'_{n+1} = w'_nv_n$. As in the
proof of the case $u = 1$ it follows that
$\relpos(P_{J_{\infty}},{}^{w'_{\infty}}P_{K_{\infty}}) = w'_{\infty}$
and therefore $y_{\infty} = w'_{\infty}$ because $z_{\infty} \in
Z^{y_{\infty}}_{J_{\infty},w_0^J}$. We have $w'_{\infty} =
w'v_0v_1\cdots v_{{\infty}-1}$ with $\delta^{-1}(v_n) \leq u_n$. If we
set $v := v_0v_1\cdots v_{{\infty}-1}$, we therefore have
$\delta^{-1}(v) \leq u$ by \eqref{addu} and this finishes the
proof.\hfill$\square$
\end{segment}

\begin{segment}\label{commutelemma}
\begin{segproc}{Lemma}
Define $\delta\colon W_J \to W_K$ as in \eqref{definedelta}. Let $J'
\subset J$, $w' \in {}^{J'}W$, $b \in B$ and $u \in W_{J'}$. We set
$K' = \delta(J') = w_0^JF(J')(w_0^J)^{-1} \subset K$. Then
there exists $v \in W_{K'}$ such that $v \leq \delta(u)$ and such that
\[
w'w_0^JF(b)F(u) \in {}^{w'}U_{K'}w'vw_0^JB.
\]
\end{segproc}

\begin{proof}
We show the claim by induction on $\ell(u)$. If $u = 1$, nothing has to
be shown. Now write $u = u_1s$ with $\ell(u_1) < \ell(u)$ and $s \in J'$. Then by
induction hypothesis there exists $v_1 \in W_{K'}$ such that
$v_1 \leq \delta(u_1)$ and such that
\[
w'w_0^JF(b)F(u) = w'w_0^JF(b)F(u_1)F(s) \in {}^{w'}U_{K'}w'v_1w_0^Jb'F(s)
\]
for some $b' \in B$. Write $b' = b_2b_1$
where $b_1 \in U_{P_{\{F(s)\}}}T$ and $b_2 \in U_{\{F(s)\}}$
\eqref{multunipotent}. Then $b'F(s) = b_2F(s)b'_1$ for some $b'_1 \in
U_{P_{\{t\}}}T \subset B$. Therefore, we can replace $b'$ by $b_2$ and
hence assume $b' \in U_{\{F(s)\}}$.

Now we have ${}^{w_0^J}U_{\{F(s)\}} = U_{\{\delta(s)\}} \subset
U_{K'}$ and therefore either ${}^{w'v_1w_0^J}U_{\{F(s)\}}$ or
${}^{w'v_1w_0^J}U^{\opp}_{\{F(s)\}}$ is contained in ${}^{w'}U_{K'}$. If
${}^{w'v_1w_0^J}U_{\{F(s)\}} \subset {}^{w'}U_{K'}$, we have
\[
{}^{w'}U_{K'}w'v_1w_0^Jb'F(s) \subset {}^{w'}U_{K'}w'v_1w_0^JF(s) =
{}^{w'}U_{K'}w'v_1\delta(s)w_0^J.
\]
In this case we set $v = v_1\delta(s)$ and have $v \leq \delta(u_1)\delta(s) =
\delta(u)$. Otherwise we note that $b'F(s) \in U^{\opp}_{F(s)}B$ and therefore
\[
{}^{w'}U_{K'}w'v_1w_0^Jb'F(s) \subset {}^{w'}U_{K'}w'v_1w_0^Jb'F(s)B.
\]
Then $v := v_1 \leq \delta(u_1) < \delta(u)$. Therefore the lemma is proved
in both cases.
\end{proof}
\end{segment}

%--------------------------------------------------------------------------

\section{Specialization of $(B \times B)$-orbits}

\begin{segment}\label{Frobliniso}
Let $\Htilde$ be any algebraic group over an algebraically closed
extension $k$ of $\FF_q$ and let $X$ and $Y$ be two varieties
with $\Htilde$-action. Let $\psi\colon X \to Y$ be a morphism of
varieties which is a homeomorphism and which satisfies
\begin{equation}\label{Froblin}
\psi(\htilde\cdot x) = F(\htilde)\psi(x)
\end{equation}
for all $x \in X$ and $\htilde \in \Htilde$. Then $\psi$ induces a
bijection of $\Htilde$-orbits which is compatible with taking
closures.

If $\Htilde$ is of the form $H \times H$ and if the condition
\eqref{Froblin} is replaced by
\begin{equation}\label{Frobsemilin}
\psi((h,h')\cdot x) = (h,F(h'))\cdot \psi(x),
\end{equation}
the analogous statement holds.
\end{segment}

\begin{segment}\label{defX}
The group $P_J^{\opp} \times P_J$ acts on the right on $G \times G
\times L_J$ by
\begin{equation}\label{actionofpar}
(g,g',z) \cdot (s,t) = (gs,g't,\pi_J^{\opp}(s)^{-1}zF(\pi_J(t))).
\end{equation}
We denote by
\[
X_J := (G \times G \times L_J)/(P_J^{\opp} \times P_J)
\]
the fppf-quotient of this action.

For any element $(g,g',z) \in G \times G \times L_K$ we denote by
$[g,g',z]$ its image in $X_J$.

The variety $X_J$ carries a left $G \times G$-action by
\[
(h,h')\cdot [g,g',z] = [hg,h'g',z].
\]
\end{segment}

\begin{segment}\label{compareZ}
\begin{segproc}{Lemma}
The map
\[
[g,g',z] \sends ({}^{g'}P_J,{}^gP_J^{\opp},[gzF(g')^{-1}])
\]
defines a $G \times G$-equivariant isomorphism
\[
\rho\colon X_J \liso Z_J
\]
of varieties with transitive $G \times G$-action.
\end{segproc}

\begin{proof}
We have
\begin{align*}
\relpos({}^{g}P_J^{\opp},{}^{gzF(g')^{-1}}F({}^{g'}P_J)) &=
\relpos(P_J^{\opp},{}^zP_J) \\
&= \relpos(P_J^{\opp},P_J) = w_0^J
\end{align*}
and hence $\rho$ lands in $Z_J$.

It is straight forward to check that $\rho$ is well-defined (use
$F(P_J) = P_J$) and $G
\times G$-equivariant. Now we show that $\rho$ is a monomorphism. As $\rho$ is
$G \times G$-equivariant, it suffices to show that $\rho([1,1,z_1]) =
\rho([g,g',z_2])$ implies $[1,1,z_1] = [g,g',z_2]$. But $\rho([1,1,z_1]) =
\rho([g,g',z_2])$ implies that $g \in P_J^{\opp}$ and $g' \in P_J$ and
therefore $[g,g',z_2] = [1,1,z_3]$ for some $z_3 \in L_J$. But
$\rho([1,1,z_1]) = \rho([1,1,z_3])$ means that $z_3 \in
U_{P_J^{\opp}}z_1U_{P_J}$. As $z_1, z_3 \in L_J$, this implies $z_1 = z_3$.

To prove the transitivity of the action of $G \times G$ on $Z_J$, it
suffices to show the following: For any $[g] \in U_{P_J^{\opp}}\backslash G /
U_{P_J}$ such that $\relpos(P_J^{\opp},{}^gP_J) = w_0^J$
there exist $(h,h') \in P_J^{\opp} \times P_J$ such that
$[hgF(h')^{-1}] = [1]$. The relation $\relpos(P_J^{\opp},{}^gP_J) =
w_0^J$ means that $P_J^{\opp}$ and ${}^gP_J$ are in opposition, that
is, their intersection is a common Levi subgroup $L$. There exists a
(unique) $u \in U_{P_J^{\opp}}$ such that ${}^uL = L_J$ and therefore
$P_J^{\opp} \cap {}^{ug}P_J = L_J$ which implies $ug \in P_J$. Now
choose $h' \in P_J$ such that $F(h') = ug$ and set $h = 1$. Then
$[hgF(h')^{-1}] = [gg^{-1}u^{-1}] = [1]$. Note that this argument
shows that the action for any affine $k$-scheme $S$ the action of
$G(S) \times G(S)$ is transitive on $Z_J(S)$.

As $G \times G$ acts transitively on $Z_J$, $\rho$ is
surjective. Moreover, as $G(k[\eps]) \times G(k[\eps])$ acts also
transitively on $Z_J(k[\eps])$ for $k[\eps] = k[T]/(T^2)$, $\rho$ is
also surjective on tangent spaces. As $X_J$ and $Z_J$ are smooth, we
therefore have a smooth surjective monomorphism and therefore an isomorphism.
\end{proof}
\end{segment}

\begin{segment}\label{Springertheorem}
For $x,w \in W$ we denote by $\Sigma^{x,w}$ the $B \times B$-orbit of
$[x,w,1]$ in $X_J$.

\begin{proclamation}{Theorem}
Every $(B \times B)$-orbit of $X_J$ is of the form $\Sigma^{x,w}$ for a
pair $(x,w) \in W \times W$. We have $\Sigma^{x,w} = \Sigma^{x',w'}$
if and only if there exists $u \in W_J$ such that
\[
x' = xF(u)^{-1}, \qquad w'u = w.
\]
For $x,x' \in W^J$ and $w,w' \in W$ the orbit $\Sigma^{x',w'}$ is
contained in the closure of the orbit of $\Sigma^{x,w}$ if and only if
there exists $u \in W_J$ such that
\begin{equation}\label{Springerorder}
xu^{-1} \leq x', \qquad F(w')u \leq F(w).
\end{equation}
\end{proclamation}

\begin{proof}
We set
\[
X^{\ad}_J := (G^{\ad} \times G^{\ad} \times
(L^{\ad}_J)^{\ad})/((P^{\opp}_J)^{\ad} \times P_J^{\ad})
\]
where the action of $(P^{\opp}_J)^{\ad} \times P_J^{\ad}$ is defined
as in \eqref{actionofpar}. Here $(L_J^{\ad})^{\ad}$ is the adjoint
group of the image of $L_J$ in $G^{\ad}$. The smooth surjective
projection morphism $X_J \to X_J^{\ad}$ induces a bijection of $(B
\times B)$-orbits in $X_J$ and $(B^{\ad} \times B^{\ad})$-orbits in
$X_J^{\ad}$ which preserves the closure of orbits. Hence it suffices
to show the analogous assertion for $(B^{\ad} \times B^{\ad})$-orbits in
$X_J^{\ad}$. In particular we can assume that $G = G^{\ad}$ and therefore
\[
X^{\ad}_J := (G \times G \times L^{\ad}_J)/(P^{\opp}_J \times P_J)
\]
with the action as in \eqref{actionofpar}.

Define a linear version of $X_J^{\ad}$ as
\[
\Gbar_J = (G \times G \times L_J^{\ad})/(P^{\opp}_J \times P_J),
\]
now with the action given by
\[
(g,g',z)\cdot (s,t) = (gs,g't,\pi^{\opp}_J(s)^{-1}z\pi_J(t)).
\]
This is the $(G \times G)$-orbit of the wonderful compactification of
de Concini and Procesi corresponding to $J$ (cf.~\cite{Sp}).

The morphism of varieties
\begin{align*}
\psi\colon X^{\ad}_J &\to \Gbar_J, \\
[g,g',z] &\sends [g,F(g'),z]
\end{align*}
is a homeomorphism such that for all $(h,h') \in G \times G$ and $x \in
X^{\ad}_J$ the equality \eqref{Frobsemilin} holds. Therefore the
theorem follows from the analogous result of Springer for the wonderful
compactification (\cite{Sp} 1.3 and 2.2).
\end{proof}
\end{segment}

\begin{segment}\label{defineBBinZ}
For $(x,w) \in W \times W$ we denote the image of $\Sigma^{x,w}$ in
$Z_J$ under the isomorphism $X_J \iso Z_J$ again by $\Sigma^{x,w}$,
i.e., as subvariety of $Z_J$ we have
\[
\Sigma^{x,w} = (B \times B)\cdot ({}^wP_J,{}^xP_J^{\opp},
[xF(w)^{-1}]).
\]
\end{segment}

\begin{segment}\label{BBG}
To relate $(B \times B)$-orbits and $G_{\rm diag}$-orbits of $Z_J$ we
need the following lemma:

\begin{proclamation}{Lemma}
For any $w \in {}^JW$ we have
\[
Z_J^w = G_{\rm diag}\cdot \Sigma^{ww_0^J,1}.
\]
\end{proclamation}

\begin{proof}
For any point $z \in Z_J$ we have
\[
G_{\rm diag}\cdot(B \times B)\cdot z = G_{\rm diag}\cdot((1)
\times B)\cdot z.
\]
Therefore it suffices to show that
\[
z(w,b) := (P_J,{}^{ww_0^J}P_J^{\opp},[ww_0^JF(b)^{-1}])
\]
is contained in $Z^w_J$ for all $b \in B$. As we have
${}^{ww_0^J}P_J^{\opp} = {}^wP_K$ and $w \in {}^JW$, this follows from
\eqref{keycorollary}.
\end{proof}
\end{segment}

\begin{segment}\label{closure}
\begin{segproc}{Lemma}
Let $w \in W$. Then
\[
\overline{\Sigma^{ww_0^J,1}} = \bigcup_{\substack{x \in W \\ x \leq
    w}} \Sigma^{xw_0^J,1}.
\]
\end{segproc}

\begin{proof}
We use the notations of the proof in \eqref{Springertheorem}. Again we
can replace $G$ by $G^{\ad}$, $X_J$ by $X_J^{\ad}$, and every $(B
\times B)$-orbit $\Sigma^{x,w}$ by its image in $X_J^{\ad}$ which we
denote again by $\Sigma^{x,w}$. The
homeomorphism $\psi$ maps the orbit $\Sigma^{ww_0^J,1}$ to the $(B
\times B)$-orbit of $[ww_0^J,1,1]$ in $\Gbar_J$.

Set $K = w_0Jw_0$. By \cite{Sp}~1.2 there exists a unique isomorphism
\[
\sigma\colon \Gbar_J \liso \Gbar_K
\]
such that
\begin{alignat*}{2}
\sigma((h,h')\cdot x) &= (h',h)\cdot\sigma(x), &&\qquad\text{for $h,h'
  \in G$, $x \in \Gbar_J$,} \\
\sigma([1,1,1]) &= [w_0,w_0,1].
\end{alignat*}
The image of the $(B \times B)$-orbit of $[ww_0^J,1,1]$ under $\sigma$ is
the $(B \times B)$-orbit of $[w_0,ww_0w_{0,J}w_0,1] = [w_0^K,w,1]$ in
$\Gbar_K$. Let $y \in W^K$ and $x \in W$ be such that the $(B \times
B)$-orbit of $[y,x,1]$ in $\Gbar_K$ is contained in the closure
of the $(B \times B)$-orbit of $\cdot[w_0^K,w,1]$. By \cite{Sp}~2.2
this means that there exists a $v \in W_K$ such that
\begin{equation}\label{intspec}
w_0^Kv^{-1} \leq y, \qquad xv \leq w.
\end{equation}
But $w_0^K$ is the maximal element in $W^K$ and therefore $y \leq
w_0^K$. Moreover, for any $v \in W_K$ we have $w_0^K \leq
w_0^Kv^{-1}$. Hence the relations \eqref{intspec} are equivalent to $v
= 1$, $y = w_0^K$, and $x \leq w$.

Now the $(B \times B)$-orbit of $[w_0^K,x,1]$ in $\Gbar_K$ corresponds
via $\sigma \circ \psi$ to $\Sigma^{xw_0^K,1}$ which shows the lemma.
\end{proof}
\end{segment}

%--------------------------------------------------------------------------

\section{A partial order on ${}^JW$}

\begin{segment}
In this chapter, we fix subsets $J$ and $K$ of $I$ and denote by
$\delta\colon W_J \to W_K$ an automorphism such that $\delta(J) =
K$. In particular we have $\delta(u) \leq \delta(u')$ if and only if
$u \leq u'$ for $u,u' \in W_J$.

We will apply the results of this chapter with the automorphism
$\delta(u) = w_0w_{0,J}F(u)w_{0,J}w_0$.
\end{segment}

\begin{segment}\label{defineSpec}
\begin{segproc}{Definition} For elements $w,w' \in W$ we write $w
  \preceq_{J,\delta} w'$ or simply $w \preceq w'$ if there exists an
  element $u \in W_J$ such that $u^{-1}w\delta(u) \leq w'$.
\end{segproc}

This order has also been examined in \cite{He} in the case that the
isomorphism $\delta$ is induced by an isomorphism $\dgtilde\colon W
\iso W$ with $\dgtilde(I) = I$. All of the following results in this
paragraph are variants of the results in loc.~cit., and the proofs are
easy modifications of the arguments there.
\end{segment}

\begin{segment}
By definition we have
\begin{equation}\label{SpecBruhat}
w \leq w' \implies w \preceq w'.
\end{equation}
If $w \in {}^JW$, we have $\ell(u^{-1}w\delta(u)) \geq \ell(u^{-1}w) -
\ell(\delta(u)) = \ell(u) + \ell(w) - \ell(u) = \ell(w)$. Therefore we
see that for $w \in {}^JW$ we have
\begin{equation}\label{Speclength}
w \preceq w' \implies \ell(w) \leq \ell(w').
\end{equation}
\end{segment}

\begin{segment}\label{Bruhatone}
Now we give three lemmas on the Bruhat order which are all proved in
\cite{He}.

\begin{proclamation}{Lemma}
Let $x,w \in W$. Then the subset
\[
\set{y \in W}{wy \leq x}
\]
of $W$ contains a smallest element $y_{\min}$ and a largest element
$y_{\max}$. Moreover we have
\[
\ell(y_{\min}) = \ell(w) - \ell(wy_{\min}), \quad \ell(y_{\max}) =
\ell(w) + \ell(wy_{\max}).
\]
\end{proclamation}
\end{segment}

\begin{segment}\label{Bruhattwo}
\begin{segproc}{Lemma}
Let $x', w, w' \in W$ such that $w \leq w'$.
\begin{assertionlist}
\item There exists $x \leq x'$ such that $xw \leq x'w'$.
\item There exists $x \leq x'$ such that $x'w \leq xw'$.
\end{assertionlist}
\end{segproc}
\end{segment}

\begin{segment}\label{Bruhatthree}
\begin{segproc}{Lemma}
Let $x \in {}^JW$, $u \in W$ such that $\ell(xu) = \ell(x) +
\ell(u)$. We write $xu = u'x'$ with $u' \in W_J$ and $x' \in
{}^JW$. Then for any $u'_1 \leq u'$ there exists a $u_1 \leq u$ such
that $xu_1 = u'_1x'$.
\end{segproc}
\end{segment}

\begin{segment}\label{Bruhatfour}
\begin{segproc}{Lemma}
Let $w \in {}^JW$, $u,v \in W_J$ such that $v
\leq u$. Then there exists $x \leq v$ and a reduced decomposition
$x = s_1\dots s_r$ such that
\begin{equation}\label{lengthsame}
\ell(s_i\dots s_1 w \delta(s_1)\dots\delta(s_i)) = \ell(w), \qquad\qquad
\text{for all $i = 1,\dots,r$}
\end{equation}
and such that
\[
x^{-1}w\delta(x) \leq u^{-1}w\delta(v).
\]
\end{segproc}

\begin{proof}
The proof is by induction on $\#J$. Assume that the statement holds
for all $J' \subset I$ with $\#J' < \#J$ and for the isomorphism
$\delta' = \delta\restricted{W_{J'}}\colon W_{J'} \iso
W_{\delta(J')}$. Then the statement is proven for $J$, $K$, and
$\delta$ by induction on $\ell(u)$.

Write $w = \wbar b$ for $\wbar \in {}^JW^K$ and $b \in W_K \cap {}^{K'}W$
where $K' = K \cap {}^{\wbar^{-1}}J$ \eqref{Howlett}. We also set $J'
= J \cap {}^{\wbar}K = {}^{\wbar}K'$. 

We first consider the case that $u \in W_{J'}$ and therefore $v \in
W_{J'}$. As $J' \subset J$, we have $w \in {}^{J'}W$. If $J'
\subsetneqq J$, we are done by induction hypothesis. If $J' = J$ and
hence $K' = K$, we have $b = 1$ and therefore $w \in {}^JW^K$. But
this implies $u^{-1}w\delta(v) \geq w$, i.e., we can choose $x = 1$.

If $u \notin W_{J'}$, we can write $u = u'u_1$ with $u_1 \in {}^{J'}W$
and $u' \in W_{J'}$ such that $\ell(u') < \ell(u)$. We also write $v =
v'v_1$ with $v_1 \leq u_1$ and $v' \leq u'$ and $\ell(v) = \ell(v') +
\ell(v_1)$. By induction hypothesis there exists $x' \leq v'$ and a
reduced decomposition $x' = s_1\dots s_{r'}$ with $\ell(s_i\dots
s_1w\delta(s_1)\dots\delta(s_i)) = \ell(w)$ for all $i = 1,\dots,r'$ such that
\[
w' := x^{\prime-1}w\delta(x') \leq u^{\prime -1}w\delta(v').
\]
Now let $x_1 \leq v_1$ be the element in $W$ such that $w'\delta(x_1)$
is the smallest element in $\set{w'\delta(y)}{y \leq v_1}$
\eqref{Bruhatone}. Then
\begin{equation}\label{leeq}
\ell(x_1^{-1}w') = \ell(w') - \ell(x_1)
\end{equation}
and
\begin{equation}\label{subeq}
w'\delta(x_1) \leq w'\delta(v_1) = u^{\prime -1}w\delta(v).
\end{equation}
Now $w\delta(v) \in \wbar W_K$ and therefore we can write $w\delta(v)
= a'\wbar b'$ with $a \in W_{J'}$ and $b' \in W_K \cap {}^{K'}W$ by
\eqref{Howlettvariant}. We now examine
\[
u^{-1}w\delta(v) = u_1^{-1}u^{\prime -1}a'\wbar b'.
\]
We have $u_1^{-1} \in W^{J'}$, $u^{\prime -1}a' \in W_{J'}$, $\wbar b'
\in {}^JW$ and $u_1^{-1}u^{\prime -1}a' \in W_J$. Therefore
\begin{align*}
\ell(u^{-1}w\delta(v))
&= \ell(u_1^{-1}u^{\prime -1}a'\wbar b') \\
&= \ell(u_1^{-1}u^{\prime -1}a') + \ell(\wbar b') \\
&= \ell(u_1^{-1}) + \ell(u^{\prime -1}a') + \ell(\wbar b') \\
&= \ell(u_1^{-1}) + \ell(u^{\prime -1}a'\wbar b') \\
&= \ell(u_1^{-1}) + \ell(u^{\prime -1}w\delta(v)).
\end{align*}
Hence it follows from \eqref{subeq} and $x_1 \leq u_1$ that
\[
(x'x_1)^{-1}w\delta(x'x_1) = x_1^{-1}w'\delta(x_1) \leq
uw\delta(v)^{-1}.
\]
Moreover, $x'x_1 \leq v$. Let $x_1 = s_1\dots s_l$ be a reduced
decomposition. For $i = 0,\dots,l$ set
\[
y_i = s_{i+1}\dots s_l (x_1^{-1}w') \delta(s_1) \dots \delta(s_i) =
s_i \dots s_1 w' \delta(s_1) \dots \delta(s_i).
\]
If we show that $\ell(y_i) = \ell(w')$ for all $i = 1,\dots,r$, we
are done. We have $\ell(y_i) \leq l +
\ell(x_1^{-1}w') = \ell(w')$ for all $i$ by \eqref{leeq}. On the other
hand, each $y_i$ is of the form $a_i^{-1}w\delta(a_i)$ for some $a_i \in
W_J$ and therefore we have $\ell(y_i) \geq \ell(a_i^{-1}w) - \ell(a_i)
= \ell(w) = \ell(w')$ as $w \in {}^JW$.
\end{proof}
\end{segment}

\begin{segment}\label{Specvariant}
\begin{segproc}{Corollary}
For $w,w' \in {}^JW$, $w \preceq w'$ if and only if there exists $u
\in W_J$ and a reduced decomposition $u = s_1\dots s_r$ such that $\ell(s_i\dots
s_1w\delta(s_1)\dots\delta(s_i)) = \ell(w)$ for all $i = 1,\dots,r$
and such that $u^{-1}w\delta(u) \leq w'$.
\end{segproc}
\end{segment}

\begin{segment}\label{SpecCoroll}
\begin{segproc}{Corollary}
Let $w,w' \in {}^JW$, $u,v \in W_J$ with $v \leq
u$. Assume that $uw'\delta(v)^{-1} \leq w$. Then $w' \preceq w$.
\end{segproc}

\begin{proof}
By \eqref{Bruhatfour} we know that there exists $x \in W_J$ such that
$x^{-1}w'\delta(x) \leq uw'\delta(v)^{-1} \leq
w$, hence $w' \preceq w$.
\end{proof}
\end{segment}

\begin{segment}\label{LemmaSpec1}
\begin{segproc}{Lemma}
Let $w,w' \in {}^JW$ such that $w
\preceq_{J,\delta} w'$. Then there exist $u,u' \in W_J$ such that
\[
uw \leq u'w'\delta(u')^{-1}\delta(u) =: w'_1
\]
and such that $w'_1 \in {}^JW$.
\end{segproc}

\begin{proof}
By definition there exists $u_1 \in W_J$ such that
$u_1^{-1}w\delta(u_1) \leq w'$. By \eqref{Bruhattwo} there exists a
$v_1 \leq \delta(u_1)^{-1}$ such that
$u_1^{-1}w\delta(u_1)\delta(u_1)^{-1} = u_1^{-1}w \leq w'v_1$. Now
$v_1 \leq \delta(u_1)^{-1}$ implies $\delta^{-1}(v_1) \leq u_1^{-1}$
and therefore $\delta^{-1}(v_1)w \leq u_1^{-1}w$ as $w \in
{}^JW$. Hence we see that there exists $v_1 \in W_K$ such that
\begin{equation}\label{Spec1eq1}
\delta^{-1}(v_1)w \leq w'v_1.
\end{equation}
Let $v_1 \in W_K$ be a minimal element such that \eqref{Spec1eq1}
holds. Then $\ell(w'v_1) = \ell(w') + \ell(v_1)$. Now
write $w'v_1 = v'w'_1$ with $v' \in W_{J}$ and $w'_1 \in
{}^{J}W$. As $v'w'_1 \geq \delta^{-1}(v_1)w$, there exists by
\eqref{Bruhattwo} an element $v'_1 \leq v'$ such that
\begin{equation}\label{Spec1eq2}
w'_1 \geq v_1^{\prime -1}\delta^{-1}(v_1)w.
\end{equation}
By \eqref{Bruhatthree} applied to $w'v_1 = v'w'_1$ and $v'_1 \leq v'$, we
see that there exists $v_2 \leq v_1$ such that $v'_1w'_1 = w'v_2$. As we
have $\ell(v'_1x') = \ell(v'_1) + \ell(x')$, \eqref{Spec1eq2} implies
$v'_1w'_1 \geq \delta^{-1}(v_1)w$. Therefore we have
\begin{equation}
w'v_2 = v'_1w'_1 \geq \delta^{-1}(v_1)w \geq \delta^{-1}(v_2)w.
\end{equation}
By the minimality of $v_1$ we have $v_1 = v_2$ and hence $v'_1 =
v'$. Therefore $w'_1 \geq v^{\prime -1}\delta^{-1}(v_1)w$. Now set $u
:= v^{\prime -1}\delta^{-1}(v_1) \in W_J$ and $u' := v^{\prime
  -1}$. Then we have $w' = v'w'_1v_1^{-1} =
u^{\prime-1}w'_1\delta(u^{-1})\delta(u')$.
\end{proof}
\end{segment}

\begin{segment}\label{lengthequal}
\begin{segproc}{Lemma}
Let $w,w' \in {}^JW$ such that $w \preceq w'$ and $\ell(w) =
\ell(w')$. Then $w = w'$.
\end{segproc}

\begin{proof}
By \eqref{Specvariant} there exist $u \in W_J$ and a reduced
decomposition $u = s_1\dots s_r$ such that $\ell(s_i\dots
s_1w\delta(s_1)\dots\delta(s_i)) = \ell(w)$ for all $i = 1,\dots,r$
and such that $u^{-1}w\delta(u) \leq w'$. Moreover we have by
\eqref{Speclength}
\[
\ell(w) \leq \ell(u^{-1}w\delta(u)) \leq \ell(w') \leq \ell(w)
\]
and hence $u^{-1}w\delta(u) = w' \in {}^JW$. We show that $w = w'$ by
induction on $\ell(u)$.

Assume that $u = s \in J$ is a simple reflection. As $\ell(w) =
\ell(sw\delta(s))$ and $\ell(sw) = \ell(w) + 1$, we have $w\delta(s)
\notin {}^JW$ and therefore $w\delta(s) = s'w$ for some $s' \in J$ by
\eqref{JWremark}. As $\ell(w) = \ell(sw\delta(s)) = \ell(ss'w) =
\ell(ss') + \ell(w)$, we have $s = s'$ and therefore $sw\delta(s) =
w$.

Now assume that $\ell(u) > 1$ and write $u = s_1u_1$ with $\ell(u_1) <
\ell(u)$. Then
\[
w' = u^{-1}w\delta(u) = u_1^{-1}s_1w\delta(s_1)\delta(u_1).
\]
We claim that $s_1w\delta(s_1) \in {}^JW$. If this is shown, we are done
by induction hypothesis and induction start.

Assume $s_1w\delta(s_1) \notin {}^JW$. If $w\delta(s_1) \notin {}^JW$,
we can again write as above $s_1w\delta(s_1) = s_1s'w$ for some $s'
\in J$. As $\ell(s_1w\delta(s_1)) = \ell(w)$, we have $s_1w\delta(s_1) =
w \in {}^JW$, contradiction. Therefore $w\delta(s_1) \in {}^JW$ and
$\ell(w\delta(s_1)) < \ell(w)$. As $w' \in {}^JW$, this implies
\begin{align*}
\ell(w) + \ell(u) &= \ell(w') + \ell(u) = \ell(uw') =
\ell(w\delta(s_1)\delta(u_1)) \\
&\leq \ell(w\delta(s_1)) + \ell(\delta(u_1)) < \ell(w) + \ell(u_1) \\
&< \ell(w) + \ell(u)
\end{align*}
and we obtain again a contradiction.
\end{proof}
\end{segment}

\begin{segment}\label{Specorder}
\begin{segproc}{Proposition}
The relation $\preceq_{J,\delta}$ is a partial order on ${}^JW$.
\end{segproc}

\begin{proof}
Let $w,w' \in {}^JW$ such that $w \preceq w'$ and $w' \preceq w$. By
\eqref{Speclength} we know that $\ell(w) = \ell(w')$ and hence $w =
w'$ by \eqref{lengthequal}. This proves the asymmetry of
$\preceq_{J,\delta}$.

Now let $w_1, w_2, w_3 \in {}^JW$ with $w_3 \preceq w_2$ and $w_2
\preceq w_1$. By \eqref{LemmaSpec1} there exist $u,u' \in W_J$ such
that $uw_3 \leq u'w_2\delta(u')^{-1}\delta(u) =: w'_2 \in {}^JW$. As
$w_2 = u^{\prime-1}w'_2\delta(u)^{-1}\delta(u') \preceq w_1$, there
exists $v \in W_J$ such that
\begin{equation*}
v^{-1}w'_2\delta(u^{-1}v) \leq w_1.
\end{equation*}
As $w'_2 \in {}^JW$, the relation $uw_3 \leq w'_2$ implies that
\begin{equation}\label{ordereq1}
v^{-1}uw_3 \leq v^{-1}w'_2.
\end{equation}
Applying \eqref{Bruhattwo} to \eqref{ordereq1}, there exists $v' \leq
u^{-1}v$ such that
\[
v^{-1}uw_3\delta(v') \leq v^{-1}w'_2\delta(u^{-1}v) \leq w_1.
\]
Then \eqref{Bruhatfour} implies $w_3 \preceq w_1$, and therefore the
relation is transitive.
\end{proof}
\end{segment}

%--------------------------------------------------------------------------

\section{The specialization order}

\begin{segment}\label{defdelta}
From now on we are back in the situation of \eqref{defZ} and we set
\[
\delta\colon W_J \liso W_K, \qquad u \sends w_0^JF(u)(w_0^J)^{-1}.
\]
\end{segment}

\begin{segment}\label{orbitclosure}
We first need the following general lemma.

\begin{proclamation}{Lemma}
Let $H$ be any algebraic group acting on a variety $Z$ and let $P
\subset H$ be an algebraic subgroup such that $H/P$ is proper. Then
for any $P$-invariant subvariety $Y \subset Z$ we have
\[
H \cdot \overline{Y} = \overline{H \cdot Y}.
\]
\end{proclamation}

\begin{proof}
Clearly we have
\[
H \cdot Y \subset H \cdot \overline{Y} \subset \overline{H \cdot Y}
\]
and therefore it suffices to show that $H \cdot \overline{Y}$ is
closed in $Z$. We denote by $\pi\colon H \times Z \to Z$ the action of $H$ on
$Z$. Define an action of $P$ on $H \times Z$ by $b\cdot (h,z) =
(hb^{-1},b\cdot z)$, and denote by $H \times^P Z$ the quotient. Then
$\pi$ induces a morphism $\pgbar\colon H \times^P Z \to Z$ which can
be written as the composition
\[
H \times^P Z \liso H/P \times Z \lto Z.
\]
Here the first morphism is the isomorphism given by $[h,z] \sends
(hP, h\cdot z)$ and the second morphism is the projection. As
$H/P$ is proper, we see that $\pgbar$ is proper. Now $\overline{Y}$ is
$P$-invariant and therefore $H \times^P \overline{Y}$ is defined, and
it is a closed subscheme of $H \times^P Z$. Therefore $\pgbar(H
\times^P \overline{Y}) = H\cdot \overline{Y}$ is closed in $Z$.
\end{proof}
\end{segment}

\begin{segment}\label{closureZBB}
\begin{segproc}{Lemma}
For $w \in {}^JW$,
\begin{equation}
\overline{Z^w_J} = \bigcup_{\substack{x \in W \\ x \leq
    w}} G \cdot\Sigma^{xw_0^J,1}.
\end{equation}
\end{segproc}

\begin{proof}
We apply \eqref{orbitclosure} to the action of $G$, embedded
diagonally in $G \times G$, on $Z_J$
and to $Y = \Sigma^{ww_0^J,1}$ which is invariant under the Borel
subgroup $B \subset G$. Then we see by \eqref{BBG}
and \eqref{closure} that
\begin{align*}
\overline{Z^w_J} &= \overline{G \cdot \Sigma^{ww_0^J,1}} \\
&= G \cdot \overline{\Sigma^{ww_0^J,1}} \\
&= \bigcup_{\substack{x \in W \\ x \leq
    w}} G \cdot\Sigma^{xw_0^J,1}.
\end{align*}
\end{proof}
\end{segment}

\begin{segment}\label{closureZ}
\begin{segproc}{Theorem}
For $w \in {}^JW$,
\[
\overline{Z^w_J} = \bigcup_{\substack{w' \in {}^JW \\ w' \preceq_{J,\delta}
    w}} Z^{w'}_J.
\]
\end{segproc}

\begin{proof}
Let $w' \in {}^JW$ with $w' \preceq w$, i.e., there exists $u \in W_J$
such that $u^{-1}w'\delta(u) \leq w$. By \eqref{closureZBB} we have
\[
G \cdot \Sigma^{u^{-1}w'\delta(u)w_0^J,1} \subset
\overline{Z^w_J}.
\]
By \eqref{BBG} we have
\begin{align*}
Z^{w'}_J &= G \cdot (P_J,{}^{w'w_0^J}P_J^{\opp}, [w'w_0^J]) \\
&= G \cdot (P_J,{}^{w'}P_K, [w'w_0^J]) \\
&= G \cdot (P_J, {}^{u^{-1}w'\delta(u)}P_K,
[u^{-1}w'w_0^JF(u)]) \\
&= G \cdot (P_J, {}^{u^{-1}w'\delta(u)}P_K, [u^{-1}w'\delta(u)w_0^J])
\end{align*}
where the third equality holds as ${}^{u^{-1}}P_J = P_J$ and
${}^{\delta(u)}P_K = P_K$. Hence $Z^{w'}_J \subset \overline{Z^w_J}$.

Conversely, let $z \in \overline{Z^w_J}$, say $z \in G \cdot
\Sigma^{xw_0^J,1}$ for some $x \in W$ with $x \leq w$
\eqref{closureZBB}. We want to show that $z \in Z_J^{w'}$ for some $w'
\in {}^JW$ with $w' \preceq w$. For this we can replace $z$ by some
element in the same $G$-orbit and hence we can assume that
$z = (P_J,{}^xP_K, [xw_0^JF(b)^{-1}])$ for some $b \in B$.

We write $x = ux'$ with $u \in W_J$ and $x' \in
{}^JW$. By \eqref{keylemma} there exists $v \in W_J$ with $v \leq u$
such that $z \in Z^{u^{-1}x\delta(v)}$. If we set $w' =
u^{-1}x\delta(v)$, we have $uw'\delta(v)^{-1} = x \leq w$ and
therefore $w' \preceq w$ by \eqref{SpecCoroll}.
\end{proof}
\end{segment}

%--------------------------------------------------------------------------
\bigskip\bigskip

\section{Applications to the Ekedahl-Oort stratification}

\begin{segment}
Let $g \geq 1$ be an integer. In this section we apply the results to
the Ekedahl-Oort stratification of the moduli space $\Acal_g$ of principally
polarized abelian varieties of dimension $g$ in characteristic
$p$. In fact all of the following results can also be applied to
arbitrary good reductions to characteristic $p > 2$ of Shimura
varieties of PEL-type (cf.~\cite{MW}, Sections~7.10-14).
\end{segment}

\begin{segment}\label{symplectic}
For this we consider the variety $Z_{G,J}$ \eqref{defZ} in a special
case: Let $\lrangle$ be the standard symplectic pairing on $V =
\FF_p^{2g}$ given by the matrix $J = \twosmallmatrix{0}{J'}{-J'}{0}$
where $J'$ is the matrix $(a_{ij})$ with $a_{ij} = \delta_{i,g+1-j}$.
Let $G = \Sp(V,\lrangle)$ be the group of symplectic isomorphisms of
$(V,\lrangle)$. Therefore the elements in $G$ are of the form
$\twosmallmatrix{A}{B}{C}{D}$ where $A$, $B$, $C$, and $D$ are $(g
\times g)$ matrices satisfying
\[
{}^tAJ'C - {}^tCJ'A = {}^tBJ'D - {}^tDJ'B = 0, \qquad {}^tAJ'D -
{}^tCJ'B = J'.
\]
The Frobenius $F\colon G \to G$ is given by $(a_{ij}) \sends
(a^p_{ij})$. As a maximal torus $T$ we choose the group of diagonal
matrices in $G$ and as the Borel subgroup $B$ we choose the subgroup
of matrices $\twosmallmatrix{A}{B}{0}{D} \in G$, where $A$ and $D$ are upper
triangular.

The normalizer of $T$ in $G$ is given by the subgroup of monomial
matrices in $GL_{2g}$ which are contained in $G$ and therefore we can
identify $W$ with the group of permutations $w \in
S_{2g}$ such that
\begin{alignat}{2}\label{symplperm}
w(i) + w(2g+1-i) &= 2g+1 & & \qquad\text{for all $i = 1,\dots,g$.}
\end{alignat}
The induced action of the Frobenius $F$ on $W$ is trivial.
An easy and elementary calculation shows that the set $I$ of
simple reflections with respect to $(T,B)$ consists of
$\{s_1,\dots,s_g\}$ with
\[
s_i =
\begin{cases} \tau_i\tau_{2g-i}, &\text{for $i = 1,\dots,g-1$;}\\
\tau_g, &\text{for $i = g$.}
\end{cases}
\]
where $\tau_j \in S_{2g}$ denotes the transposition of $j$ and $j+1$.
The longest element $w_0 \in W$ is given by the permutation $i \sends
2g+1-i$. It follows from \eqref{symplperm} that $w_0$ lies in the
center of $W$. In particular, $K = {}^{w_0}J = J$.

Let $J = \{s_1,\dots,s_{g-1}\} \subset I$. Then the standard parabolic
subgroup $P_J$ consists of the matrices $\twosmallmatrix{A}{B}{C}{D}$
in $G$ with $C = 0$, and $W_J$ is the subgroup of those permutation $w
\in W$ such that $w(\{1,\dots,g\}) = \{1,\dots,g\}$. The map
\[
W_J \to S_g. \qquad w \mapsto w\restricted{\{1,\dots,g\}}
\]
is a group isomorphism. The maximal element $w_{0,J}$ in $W_J$
corresponds via this map to the permutation $i \sends g+1-i$ in $S_g$.
The map $\delta$ from \eqref{defdelta} has therefore in this case the
form
\begin{equation}\label{deltasympl}
\delta\colon W_J \iso W_J, \quad w \sends w_{0,J}ww_{0,J}.
\end{equation}

The set ${}^JW$ consists of those elements $w \in W$ such that
\[
w^{-1}(1) < w^{-1}(2) < \dots < w^{-1}(g).
\]
Of course, this implies $w^{-1}(g+1) < \dots < w^{-1}(2g)$. For
two permutations $w$ and $w'$ in ${}^JW$ we have $w \leq w'$ if and
only if $w^{-1}(i) \leq w^{\prime-1}(i)$ for all $i = 1,\dots,g$.

If $\Sigma = \{j_1 < \dots < j_g\} \subset \{1,\dots,2g\}$ is a subset of
$g$ elements such that either $i \in \Sigma$ or $2g+1-i \in \Sigma$
for all $i = 1,\dots,g$, we get a corresponding element $w_{\Sigma}
\in {}^JW$ by setting $w^{-1}(i) = j_i$.
The sets of these $\Sigma$'s is in bijection with $\{0,1\}^g$ by
associating to $\Sigma$ the tuple $(\epsilon_1,\dots,\epsilon_g)$ with
\[
\epsilon_i =
\begin{cases}
0, &\text{if $i \in \Sigma$;}\\
1, &\text{otherwise.}
\end{cases}
\]
The length of such an element $(\epsilon_1,\dots,\epsilon_g)$ is equal
to
\[
\sum_{i=1}^g \epsilon_i(g+1-i).
\]
\end{segment}

\begin{segment}
The moduli space $\Acal_g$ is a Deligne-Mumford stack. Let
$(\Xcal,\lambda)$ be the universal abelian (relative) scheme over
$\Acal_g$ and $f\colon \Xcal \to \Acal_g$ its structure morphism.
The first de Rham cohomology $M = H^1_{\DR}(\Xcal/\Acal_g)$
is endowed with the structure of an $F$-Zip with symplectic
structure (cf.\ \cite{MW},~\S 7): $M$ is a locally free
$\Ocal_{\Acal_g}$-module of rank $2g$ and the principal polarization
induces a symplectic pairing $\beta$ on $M$. There are two canonical
locally direct summands of rank $g$, namely
\[
C^1 := f_*\Omega^1_{\Xcal/\Acal_g}, \qquad D_0 :=
R^1f_*(\Hscr(\Omega^{\bullet}_{\Xcal/\Acal_g})).
\]
The Cartier isomorphism induces $\Ocal_{\Acal_g}$-linear isomorphisms
\begin{align*}
\varphi_0\colon (M/C^1)^{(p)} \cong (R^1f_*\Ocal_{\Xcal})^{(p)} &\liso
D_0 \\
\varphi_1\colon (C^1)^{(p)} &\liso
f_*(\Hscr^1(\Omega^{\bullet}_{\Xcal/\Acal_g})) \cong M/D_0.
\end{align*}
Moreover, $C^1$ and $D_0$ are both totally isotropic with respect to
the symplectic pairing $\beta$.

We obtain a morphism
\[
\pi\colon \Acal_g \lto [G\backslash Z_{G,J}]
\]
where $G$ and $J$ as defined in \eqref{symplectic}. Here the right
hand side is the quotient in the sense of stacks. As the $G$-orbits of
$Z_{G,J}$ are parametrized by elements in ${}^JW$, the underlying
topological space of $[G\backslash Z_{G,J}]$ is in bijection with
${}^JW$ and for $w \in {}^JW$ we denote the corresponding point in
$[G\backslash Z_{G,J}]$ by $C(w)$. It is the image of the
corresponding orbit $Z_{G,J}^w$. Moreover, a point $C(w)$ is
contained in the closure of $\{C(w')\}$ if and only if $Z_{G,J}^w
\subset \overline{Z_{G,J}^{w'}}$ (see e.g.~\cite{WdOortstrata}~(4.4)).
Each $\{C(w)\}$ is locally closed in the underlying topological space
of $[G\backslash Z_{G,J}]$ and there is a canonical structure of a
locally closed substack on $\{C(w)\}$ \cite{MW}, Section~5.6. Indeed,
this is nothing but the unique structure of a reduced substack on
$\{C(w)\}$. We denote this substack again by $C(w)$. As all $w \in
{}^JW$ are already defined over $\FF_p$, these substacks are also
defined over $\FF_p$.

For all $w \in {}^JW$ we define a locally closed substack $\Acal^w_g$
of $\Acal_g$ by the cartesian diagram
\[\xymatrix{
\Acal^w_g \ar@{^{(}->}[d] \ar[r] \ar@{}[dr]|{\square} & C(w)
\ar@{^{(}->}[d] \\
\Acal_g \ar[r]^-{\pi} & [G\backslash Z_{G,J}].
}\]

Now let $k$ be an algebraically closed field and let $x\colon \Spec(k)
\to \Acal_g$ be a $k$-valued point of $\Acal_g$. Then $x$ corresponds
to a principally polarized abelian variety $(X,\lambda)$ of dimension $g$ over
$k$. Its first de Rham cohomology $M = H^1_{\DR}(X/k)$ carries the
structure $(C^1,D_0,\varphi_0,\varphi_1,\beta)$ of a
symplectic $F$-zip as above. On the other hand, the Dieudonn\'e module $(M',F,V)$
associated to the $p$-torsion $X[p]$ of $X$ carries a symplectic
pairing $\beta'$ and we can associate a symplectic $F$-zip to
$(M',F,V,\beta')$ as follows: We set $C^{\prime1} = VM' = \Ker(F)$ and
$D'_0 = FM' = \Ker(V)$ and we denote by $\varphi'_0\colon
(M'/VM')^{(p)} \iso FM'$ and $\varphi'_1\colon (VM')^{(p)} \iso
M'/FM'$ the $k$-linear isomorphisms induced by the Frobenius linear
maps $F$ and $V^{-1}$, respectively. Then these two symplectic $F$-zips
$(M,C^1,D_0,\varphi_0,\varphi_1,\beta)(x)$ and
$(M',C^{\prime1},D'_0,\varphi'_0,\varphi'_1,\beta')(x)$ associated to
$x$ are canonically isomorphic.

If $x,x'\colon \Spec(k) \to \Acal_g$ are two $k$-valued points,
corresponding to principally polarized abelian varieties $(X,\lambda)$
and $(X',\lambda')$, respectively, then
they both factorize through the same locally closed substack
$\Acal_g^w$ if and only if the symplectic $F$-zip structures induced
on the first de Rham cohomology are isomorphic. By the above this is
equivalent to the fact that the principally quasi-polarized
Dieudonn\'e module of their $p$-torsions are isomorphic. Hence we see
that the $\Acal_g^w$ for $w \in {}^JW$ are nothing but the
Ekedahl-Oort strata, defined in \cite{Oort}.
\end{segment}

\begin{segment}
Moreover, by \cite{Wdflat} (for $p > 2$ this follows also from
\cite{WdOortstrata}) we have:

\begin{proclamation}{Theorem}
The morphism
\[
\pi\colon \Acal_g \lto [G\backslash Z_{G,J}]
\]
is faithfully flat.
\end{proclamation}
\end{segment}

\begin{segment}
In particular the theorem implies that for $w,w' \in {}^JW$, we have
$\Acal_g^w \subset \overline{\Acal_g^{w'}}$ if and only if $Z^w_{G,J} \subset
\overline{Z^{w'}_{G,J}}$. Therefore \eqref{closureZ} implies:

\begin{proclamation}{Corollary}
The following two assertions are equivalent:
\begin{assertionlist}
\item $\Acal_g^w \subset \overline{\Acal_g^{w'}}$.
\item There exists $u \in W_J$ such that $u^{-1}w\delta(u) \leq w'$
  with $\delta(u) = w_{0,J}uw_{0,J}$.
\end{assertionlist}
\end{proclamation}
\end{segment}

%==========================================================================

\end{document}